\documentclass[12pt,twoside]{amsart}
\usepackage{amsmath, amsthm, amscd, amsfonts, amssymb, graphicx}
\usepackage{enumerate}
\usepackage[colorlinks=true,
linkcolor=blue,
urlcolor=cyan,
citecolor=red]{hyperref}
\usepackage{mathrsfs}

\newtheorem{theorem}{Theorem}[section]

\newtheorem{corollary}{Corollary}[section]

\newtheorem{proposition}{Proposition}[section]
\numberwithin{equation}{section}

\input{tcilatex}

\begin{document}
\title{An $L_{p}$ norm inequality related to extremal polynomials}
\author{Abdelhamid Rehouma }
\address{Department of Mathematics, Faculty of Exact Sciences, University
Hama Lakhdar, Eloued, Algeria.}
\email{rehoumaths@gmail.com}
\author{ Herry Pripawanto Suryawan}
\address{Department of Mathematics, Faculty of Science and Technology,
Sanata Dharma University, P.O.Box 55281, Yogyakarta, Indonesia.}
\email{E-mail address: herrypribs@usd.ac.id}
\maketitle

\begin{abstract}
Let $E$ be a Jordan rectifiable curve in the complex plane and let $G$ be
the bounded component of $\mathbb{C}\backslash E.$ Now let $n\in \mathbb{N}$%
, and let $m_{n,E}$ denote the extremal constants defined by 
\begin{equation*}
m_{n,E}=\inf \left\{ \left\Vert \dfrac{D_{E,\rho }\left( z\right) }{%
D_{E,\rho }\left( 0\right) }-P_{n}\left( z\right) \right\Vert _{L^{p}\left(
G,\rho \right) }:P_{n}\left( \xi \right) =1\right\}
\end{equation*}%
where $\xi $ is a fixed complex number.where $\rho $ is a weight function, $%
D_{E,\rho }\left( \cdot \right) $ is the so called {Szeg\"{o}} function, $%
z\in G$, $p\geq 2.$ The infimum is taken over all polynomials $P_{n}$ of
degree $n$. The $L_{p}$ associated extremal polynomials $\left\{
Q_{n}\right\} _{n=1,2....}$ satisfies%
\begin{equation*}
m_{n,E}=\left\Vert \dfrac{D_{E,\rho }\left( z\right) }{D_{E,\rho }\left(
0\right) }-Q_{n}\left( z\right) \right\Vert _{L^{p}\left( G,\rho \right) }.
\end{equation*}%
We define the functions, if $p\in $ $%
\mathbb{N}
$%
\begin{equation*}
J_{n}\left( z\right) =\int_{\xi _{G}}^{z}Q_{n}^{p}\left( t\right) dt;\;z\in G
\end{equation*}%
which are of course well defined polynomials for any $n\in 
\mathbb{N}
.$Following the same convention , we define the function%
\begin{equation*}
\Phi \left( z\right) =\int\limits_{\xi _{G}}^{z}\left( \dfrac{D_{E,\rho
}\left( t\right) }{D_{E,\rho }\left( 0\right) }\right) ^{p}dt,
\end{equation*}%
Our main target in this paper is to show that when $m_{n,E}\longrightarrow
0, $ then 
\begin{equation*}
J_{n}\left( z\right) \text{ }\longrightarrow \Phi \left( z\right)
\end{equation*}%
uniformly on compact subsets of $G.$
\end{abstract}

\pagestyle{myheadings} \markboth{\centerline {}}
{\centerline {}} \bigskip \bigskip 

\section{Introduction}

Let $E$ be a {Jordan} rectifiable curve in the complex plane $%
\mathbb{C}
$. It is well known that $\mathbb{C}\backslash E$ is composed of two
connected domains, one of which is bounded. The bounded component of $%
\mathbb{C}\backslash E$ will be denoted by $G_E$ (or $G$ if only one curve
is considered).

Let $\sigma $ be a finite positive measure on the Borel $\sigma $-algebra of 
$%
\mathbb{C}
$ that is concentrated on the curve $E.$

It is well known that $\sigma $ is absolutely continuous with respect to arc
length measure $\left\vert dt\right\vert $ on $E$. Indeed, we have \cite{SZ} 
\begin{equation}
d\sigma =\rho \left( t\right) \left\vert dt\right\vert ,\ \ \ \ \
\int\limits_{E}\rho \left( t\right) \left\vert dt\right\vert <\infty
\label{eq_abs_cont}
\end{equation}%
and 
\begin{equation}
\int\limits_{E}\left\vert t\right\vert ^{n}\rho \left( t\right) \left\vert
dt\right\vert <\infty ;n\in \mathbb{N},  \label{eq_t^n}
\end{equation}%
for some non-negative function $\rho $ on $E$. This function $\rho $ will be
referred to as a weight.Moreover \cite{SL,SZ} 
\begin{equation}
\int\limits_{E}\log \left( \rho \left( t\right) \right) \left\vert \varphi
^{\prime }\left( t\right) \right\vert \left\vert dt\right\vert >-\infty
\label{eq_log}
\end{equation}%
where $\varphi $ denotes the conformal mapping of $G$ onto the disk $%
D_{R}=\left\{ w\in \mathbb{C},\left\vert w\right\vert <R\right\} $ such that 
$\varphi \left( \xi \right) =0,\varphi ^{\prime }\left( \xi \right) =1,$ for
some fixed complex number $\xi $.In these results, $E$ is still a Jordan
curve, as above, $G$ is the bounded component of $\mathbb{C}\backslash E$, $%
\varphi $ is the conformal mapping that maps $G$ onto the unit disk, with $%
|\varphi (t)|=1$ for $t\in E$, and $\psi =\varphi ^{-1}.$ We notice that
such mapping exists by \cite{JU,SL,SZ}.

For a given Jordan curve $E$, and a given weight $\rho $ as above, the
weighted {Smirnov} spaces is denoted by $L^{p}\left( G,\rho \right)
,0<p<\infty $, and it is defined by 
\begin{equation*}
L^{p}(G,\rho )=\left\{ f:G\cup E\rightarrow \mathbb{C},f\;{\text{is\
analytic\ in}}\;G,\;f\;{\text{is\ continuous\ on}}\;E,\;{\text{and}}\;\Vert
f\Vert _{L^{p}(G,\rho )}<\infty \right\} ,
\end{equation*}%
where 
\begin{equation}
\left\Vert f\right\Vert _{L^{p}\left( G,\rho \right) }=\left(
\int\limits_{E}\left\vert f\left( t\right) \right\vert ^{p}\rho \left(
t\right) \left\vert dt\right\vert \right) ^{\frac{1}{p}}.  \label{Norm}
\end{equation}%
We refer the reader to \cite{JU,SL,SZ} for further details on Smirnov spaces
and their properties.We refer the reader to \cite{JU,SL,SZ} for further
details on Smirnov spaces and their properties.(\textbf{Ahlfors-Smirnov}
regular domains\textbf{\ } ).They are defined by the condition that there
exists a constant $K>0$ such that

\begin{equation}
\left\vert E\cap D_{r}\right\vert \leq C_{r}\text{ \ \ \ \ \ }\forall
D_{r}=\left\{ z,\left\vert z\right\vert =r\right\}  \label{Condition}
\end{equation}%
where $\left\vert dt\right\vert $ denotes the arc length measure on $%
E=\partial G.$

In the following few lines, we present the {Szeg\"{o} }function, see \cite%
{SB,LB,SM,AA}. As before, we have a Jordan curve $E$ with a measure $\sigma $
concentrated on the curve $E.$ So, we have a weight $\rho $ satisfying %
\eqref{eq_abs_cont}, \eqref{eq_t^n} and \eqref{eq_log}. Then, one can
construct the so-called {Szeg\"{o}} \ function $D_{E,\rho }\left( z\right) $
associated with the curve $E$ and the weight function $\rho \left( t\right) $
with the following properties, where $G$ is the bounded component of $%
\mathbb{C}\backslash E$.

\begin{enumerate}
\item[i)] $D_{E,\rho }\left( z\right) $ is analytic in $G,$ $D_{E,\rho
}\left( z\right) $ $\neq 0$ in $G$.

\item[ii)] $D_{E,\rho }\left( 0\right) >0.$

\item[iii)] $D_{E,\rho }\left( z\right) $ has a limit \ boundary value
satisfying

\begin{equation*}
\left\vert D_{E,\rho }\left( t\right) \right\vert ^{p}\left\vert \varphi
^{\prime }\left( t\right) \right\vert =\rho \left( t\right) ,~\ t\in E~\ 
\end{equation*}%
In fact, $D_{E,\rho }$ can be defined as follows. If $\varphi $ is the
conformal mapping associated with $G$, that is the mapping that maps $G$
onto $D_{1}$, then one can define $D_{E,\rho }\left( z\right) $ $%
=D_{G}\left( \varphi \left( z\right) \right) ,$ where 
\begin{equation}
D_{G}\left( w\right) =\exp \left\{ \frac{1}{2p\pi }\int_{0}^{2\pi }\frac{%
w+e^{i\theta }}{w-e^{i\theta }}\log \frac{\rho \left( t\right) }{\left\vert
\varphi ^{\prime }\left( t\right) \right\vert }\left\vert \varphi ^{\prime
}\left( t\right) \right\vert \left\vert dt\right\vert \right\}
\label{eq_exp}
\end{equation}
\end{enumerate}

\section{Extremal polynomials}

For a Jordan curve $E$, let $D_{E,\rho }$ be the Szeg\"{o} function, as
discussed earlier. We are interested in the following well known extremal
problem 
\begin{equation*}
\mu _{p,E}=\left\Vert \dfrac{D_{E,\rho }\left( z\right) }{D_{E,\rho }\left(
0\right) }\right\Vert _{L^{p}\left( G,\rho \right) }=\inf \left\{ \left\Vert
f\right\Vert _{L^{p}\left( E,\rho \right) }:f\left( \xi _{G}\right)
=1\right\}
\end{equation*}%
where the infimum is taken over the holomorphic functions $f$ \ in the
weighted {Smirnov} spaces $L^{p}\left( G,\rho \right) .$

Now let $n\in \mathbb{N}$, and let $m_{n,E}$ denote the extremal constants
associated with the measure $\sigma $ and $G.$ That is, $m_{n,E}$ is defined
by 
\begin{equation}
m_{n,E}=\inf \left\{ \left\Vert \dfrac{D_{E,\rho }\left( z\right) }{%
D_{E,\rho }\left( 0\right) }-P_{n}\left( z\right) \right\Vert _{L^{p}\left(
G,\rho \right) }:P_{n}\left( \xi \right) =1\right\}  \label{eq_m_n}
\end{equation}%
where $\xi $ is a fixed complex number. The infimum is taken over all
polynomials $P_{n}$ of degree $n$. The $L_{p}$ extremal polynomials $\left\{
Q_{n}\right\} _{n=1,2....}$ associated with the measure $\sigma $ and the
support $G$ are defined as solutions to the best approximation problem %
\eqref{eq_m_n} such that $Q_{n}\left( \xi \right) =1,$ and%
\begin{equation*}
m_{n,E}\left( \sigma \right) =\left\Vert \dfrac{D_{E,\rho }\left( z\right) }{%
D_{E,\rho }\left( 0\right) }-Q_{n}\left( z\right) \right\Vert _{L^{p}\left(
G,\rho \right) }.
\end{equation*}%
We define the functions%
\begin{equation}
J_{n}\left( z\right) =\int_{\xi _{G}}^{z}Q_{n}^{p}\left( t\right) dt;\;z\in G
\label{eq_J_n}
\end{equation}%
which are of course well defined polynomials for any $n\in 
\mathbb{N}
$, if $p\in $ $%
\mathbb{N}
.$Following the same convention as for \eqref{eq_J_n}, we define the function%
\begin{equation}
\Phi \left( z\right) =\int\limits_{\xi _{G}}^{z}\left( \dfrac{D_{E,\rho
}\left( t\right) }{D_{E,\rho }\left( 0\right) }\right) ^{p}dt,
\label{eq_phi}
\end{equation}%
where the holomorphic function $\dfrac{D_{E,\rho }\left( z\right) }{%
D_{E,\rho }\left( 0\right) }$ satsisfies :

\begin{equation}
\left\Vert \dfrac{D_{E,\rho }\left( z\right) }{D_{E,\rho }\left( 0\right) }%
\right\Vert _{L^{p}\left( G,\rho \right) }=\inf \left\{ \left\Vert
f\right\Vert _{L^{p}\left( E,\rho \right) }:f\left( \xi \right) =1\right\}
\label{Extrpbm}
\end{equation}%
Our main target in this paper is to show that when $m_{n,E}\longrightarrow
0, $ then 
\begin{equation*}
\int_{\xi }^{z}Q_{n}^{p}\left( t\right) dt\text{ }\longrightarrow
\int\limits_{\xi }^{z}\left( \dfrac{D_{E,\rho }\left( t\right) }{D_{E,\rho
}\left( 0\right) }\right) ^{p}dt
\end{equation*}%
(Equivalently, $J_{n}(z)\rightarrow \Phi (z)$) uniformly (in $z$) on compact
subsets of $G.$ Our interest in $Q_{n}$ and $J_{n}$ is explained by the fact
that one can produce estimates of the convergence rates for these
polynomials.

\begin{theorem}
Let $E$ be a Jordan curve, and let $\sigma $ be concentrated on the curve $%
E. $ Suppose that the absolutely continuous part $\rho (t)\left\vert
dt\right\vert ,t\in E$, of $\sigma $ verifies the {Szeg\"{o} }conditions %
\eqref{eq_abs_cont} and \eqref{eq_log}. If $p\geq 2$, then%
\begin{equation}
\left\Vert J_{n}-\Phi \right\Vert _{\infty }\leq \frac{m_{n,E}}{2}\left( 
\frac{\gamma _{p,q}}{D_{E,\rho }(0)}+\left\Vert Q_{n}\right\Vert
_{L^{p}\left( G,\nu \right) }\right) ^{p-1}  \label{Equation}
\end{equation}%
where%
\begin{equation*}
\nu \left( \xi \right) \left\vert d\xi \right\vert =\rho ^{1-q}\left( \xi
\right) \left\vert d\xi \right\vert \text{ \ \ \ \ \ \ \ \ \ \ \ \ \ \ \ \ }%
\xi \in \text{\ }E
\end{equation*}%
and%
\begin{equation*}
\gamma _{p,q}=\left( \int\limits_{E}\left\vert D_{E,\rho }\left( \xi \right)
\right\vert ~^{p\left( 2-q\right) }\left\vert \varphi ^{\prime }\left( \xi
\right) \right\vert ^{1-q}\left\vert d\xi \right\vert \right) ^{\frac{1}{p}}
\end{equation*}%
and $\frac{1}{p}+\frac{1}{q}=1.$ In particular, if $m_{n,E}\left( \sigma
\right) \longrightarrow 0,$ as $n\longrightarrow \infty $, then%
\begin{equation}
\left\Vert J_{n}-\Phi \right\Vert _{\infty }\longrightarrow 0
\label{eq_prov_2}
\end{equation}%
uniformly on compact subsets of $G.$
\end{theorem}

\begin{proof}
Thanks to \eqref{eq_J_n} and \eqref{eq_phi}, we have 
\begin{align*}
\left\vert J_{n}\left( z\right) -\Phi \left( z\right) \right\vert &
=\left\vert \int\limits_{\xi }^{z}\left\{ \left( Q_{n}\left( t\right)
\right) ^{p}-\left( \frac{D_{E,\rho }(t)}{D_{E,\rho }(0)}\right)
^{p}\right\} dt\right\vert \\
& =\left\vert \int\limits_{0}^{\varphi \left( z\right) }\left\{
Q_{n}^{p}\left( \psi \left( u\right) \right) -\left( \frac{D_{E,\rho }(\psi
\left( u\right) )}{D_{E,\rho }(0)}\right) ^{p}\right\} \psi ^{\prime }\left(
u\right) du\right\vert \\
& \leq \int\limits_{0}^{\varphi \left( z\right) }\left\vert Q_{n}^{p}\left(
\psi \left( u\right) \right) -\left( \frac{D_{E,\rho }(\psi \left( u\right) )%
}{D_{E,\rho }(0)}\right) ^{p}\right\vert \left\vert \psi ^{\prime }\left(
u\right) \right\vert \left\vert du\right\vert
\end{align*}%
where $z\in G$. But $\varphi $ is the conformal mapping that maps $G$ onto
the unit disk, with $|\varphi (t)|=1$ for $t\in E$, and $\psi =\varphi
^{-1}. $The integration is carried over the segment $\left[ 0,\varphi \left(
z\right) \right] $ in $D_{1}$. Since $E$ is rectifiable, the function under
the latter integral belongs to the Hardy class $H^{1}\left( D_{1}\right) $.
Hence we obtain by the {Fejer-Riesz} inequality (cf. [$6$, Theorem $3.13$])
that%
\begin{equation*}
\left\vert J_{n}\left( z\right) -\Phi \left( z\right) \right\vert \leq \frac{%
1}{2}\int\limits_{\left\vert u\right\vert =1}\left\vert Q_{n}^{p}\left( \psi
\left( u\right) \right) -\left( \frac{D_{E,\rho }(\psi \left( u\right) )}{%
D_{E,\rho }(0)}\right) ^{p}\right\vert \left\vert \psi ^{\prime }\left(
u\right) \right\vert \left\vert du\right\vert .
\end{equation*}%
Direct computations lead to%
\begin{equation*}
\left\vert J_{n}\left( z\right) -\Phi \left( z\right) \right\vert \leq \frac{%
1}{2}\int\limits_{E}\left\vert \left( Q_{n}\left( t\right) \right)
^{p}-\left( \frac{D_{E,\rho }\left( t\right) )}{D_{E,\rho }(0)}\right)
^{p}\right\vert \left\vert dt\right\vert .
\end{equation*}%
Applying {Holder}'s inequality for $p\geq 2,$ we have%
\begin{align*}
\left\vert J_{n}\left( z\right) -\Phi \left( z\right) \right\vert & \leq 
\frac{1}{2}\int\limits_{E}\left\vert Q_{n}\left( t\right) -\frac{D_{E,\rho
}(t)}{D_{E,\rho }(0)}\right\vert \\
& \times \left\vert \left\vert D_{E,\rho }\left( t\right) \right\vert
~^{-p}\left\vert \varphi ^{\prime }\left( t\right) \right\vert
^{-1}\sum\limits_{k=0}^{p-1}\left( \frac{D_{E,\rho }(t)}{D_{E,\rho }(0)}%
\right) ^{k}\left[ Q_{n}\left( t\right) \right] ^{p-k-1}\right\vert \rho
\left( t\right) \left\vert dt\right\vert \\
& \leq \frac{1}{2}\left\Vert Q_{n}\left( t\right) -\frac{D_{E,\rho }(t)}{%
D_{E,\rho }(0)}\right\Vert _{L^{p}\left( G,\rho \right) } \\
& \times \left\Vert \left\vert D_{E,\rho }\left( t\right) \right\vert
~^{-p}\left\vert \varphi ^{\prime }\left( t\right) \right\vert
^{-1}\sum_{k=0}^{p-1}\left( \frac{D_{E,\rho }(t)}{D_{E,\rho }(0)}\right) ^{k}%
\left[ Q_{n}\left( t\right) \right] ^{p-k-1}\right\Vert _{L^{q}\left( G,\rho
\right) }
\end{align*}

where $q=\dfrac{p}{p-1}$. Now we have 
\begin{align*}
\left\vert \sum_{k=0}^{p-1}\left( \frac{D_{E,\rho }(t)}{D_{E,\rho }(0)}%
\right) ^{k}\left[ Q_{n}\left( t\right) \right] ^{p-k-1}\right\vert & \leq
\sum_{k=0}^{p-1}\left\vert \frac{D_{E,\rho }(t)}{D_{E,\rho }(0)}\right\vert
^{k}\left\vert Q_{n}\left( t\right) \right\vert ^{p-k-1} \\
& \leq \left( \left\vert \frac{D_{E,\rho }(t)}{D_{E,\rho }(0)}\right\vert
+\left\vert Q_{n}\left( t\right) \right\vert \right) ^{p-1}.
\end{align*}%
It follows that 
\begin{align*}
& \left\Vert \left\vert D_{E,\rho }\left( t\right) \right\vert
~^{-p}\left\vert \varphi ^{\prime }\left( t\right) \right\vert
^{-1}\sum_{k=0}^{p-1}\left[ \frac{D_{E,\rho }(t)}{D_{E,\rho }(0)}\right] ^{k}%
\left[ Q_{n}\left( t\right) \right] ^{p-k-1}\right\Vert _{L^{q}\left( G,\rho
\right) } \\
& \leq \left[ \int\limits_{E}\left\vert D_{E,\rho }\left( t\right)
\right\vert ~^{-qp}\left\vert \varphi ^{\prime }\left( t\right) \right\vert
^{-q}\left( \left\vert \frac{D_{E,\rho }(t)}{D_{E,\rho }(0)}\right\vert
+\left\vert Q_{n}\left( t\right) \right\vert \right) ^{p}\rho \left(
t\right) \left\vert dt\right\vert \right] ^{\frac{p-1}{p}} \\
& \leq \left( \left\Vert \frac{D_{E,\rho }(.)}{D_{E,\rho }(0)}\right\Vert
_{L^{p}\left( G,\nu \right) }+\left\Vert Q_{n}\right\Vert _{L^{p}\left(
G,\nu \right) }\right) ^{p-1},
\end{align*}%
by {Minkowski}'s inequality. Thus, we have shown that, for each $z\in G,$ 
\begin{equation*}
\left\vert J_{n}\left( z\right) -\Phi \left( z\right) \right\vert \leq \frac{%
1}{2}\left\Vert Q_{n}-\frac{D_{E,\rho }(.)}{D_{E,\rho }(0)}\right\Vert
_{L^{p}\left( G,\rho \right) }\left( \left\Vert \frac{D_{E,\rho }(.)}{%
D_{E,\rho }(0)}\right\Vert _{L^{p}\left( G,\nu \right) }+\left\Vert
Q_{n}\right\Vert _{L^{p}\left( G,\nu \right) }\right) ^{p-1}.
\end{equation*}

Since,%
\begin{equation*}
\gamma _{p,q}=\left\Vert \frac{D_{E,\rho }(.)}{D_{E,\rho }(0)}\right\Vert
_{L^{p}\left( G,\nu \right) }=D_{E,\rho }^{-1}(0)\left(
\int\limits_{E}\left\vert D_{E,\rho }\left( \xi \right) \right\vert
~^{p\left( 2-q\right) }\left\vert \varphi ^{\prime }\left( \xi \right)
\right\vert ^{1-q}\left\vert d\xi \right\vert \right) ^{\frac{1}{p}}
\end{equation*}%
and%
\begin{equation*}
m_{n,E}=\left\Vert Q_{n,p,E}\left( t\right) -\frac{D_{E,\rho }(.)}{D_{E,\rho
}(0)}\right\Vert _{L^{p}\left( G,\rho \right) },
\end{equation*}%
we have 
\begin{equation*}
\left\Vert J_{n}-\Phi \right\Vert _{\infty }\leq \frac{m_{n,E}}{2}\left( 
\frac{\gamma _{p,q}}{D_{E,\rho }(0)}+\left\Vert Q_{n}\right\Vert
_{L^{p}\left( G,\nu \right) }\right) ^{p-1},
\end{equation*}%
which proves \eqref{Equation}. Then \eqref{eq_prov_2} follows immediately.
This completes the proof.
\end{proof}

\begin{proposition}
For a weight function \ $\rho (t)$ and a contour $E,$ let 
\begin{equation}
\gamma =\sqrt{\int\limits_{E}\frac{\left\vert \varphi ^{\prime }\left(
t\right) \right\vert ^{2}}{\rho \left( t\right) }\left\vert dt\right\vert }.
\label{Longgam}
\end{equation}%
If $\gamma <\infty $, then for any function $Q\left( z\right) \in
L^{2}\left( E,\rho \right) $%
\begin{equation}
\int\limits_{\left\vert u\right\vert =1}\left\vert Q\left( \psi \left(
u\right) \right) \right\vert \left\vert du\right\vert \leq \gamma \left\Vert
Q\right\Vert _{L^{2}\left( E,\rho \right) }  \label{Inequone}
\end{equation}%
where%
\begin{equation*}
\left\Vert Q\right\Vert _{L^{2}\left( E,\rho \right) }=\sqrt{%
\int\limits_{E}\left\vert Q\left( t\right) \right\vert ^{2}\rho \left(
t\right) \left\vert dt\right\vert .}
\end{equation*}
\end{proposition}

\begin{proof}
Since the function $Q\left( z\right) $ is defined on the contour $E$, the
function $Q\left( \psi \left( u\right) \right) $ is defined when $\left\vert
u\right\vert =1.$ We then have 
\begin{equation*}
\int\limits_{\left\vert u\right\vert =1}\left\vert Q\left( \psi \left(
u\right) \right) \right\vert \left\vert du\right\vert
=\int\limits_{\left\vert u\right\vert =1}\left\vert Q\left( \psi \left(
u\right) \right) \right\vert \frac{\sqrt{\rho \left( \psi \left( u\right)
\right) \left\vert \psi ^{\prime }\left( u\right) \right\vert }}{\sqrt{\rho
\left( \psi \left( u\right) \right) \left\vert \psi ^{\prime }\left(
u\right) \right\vert }}\left\vert du\right\vert
\end{equation*}

\begin{equation*}
\leq \sqrt{\int\limits_{\left\vert u\right\vert =1}\frac{\left\vert
du\right\vert }{\rho \left( \psi \left( u\right) \right) \left\vert \psi
^{\prime }\left( u\right) \right\vert }}\sqrt{\int\limits_{\left\vert
u\right\vert =1}\left\vert Q\left( \psi \left( u\right) \right) \right\vert
^{2}\rho \left( \psi \left( u\right) \right) \left\vert \psi ^{\prime
}\left( u\right) \right\vert \left\vert du\right\vert },
\end{equation*}%
where we have used the Cauchy-Schwarz inequality to obtain the last
inequality. Hence%
\begin{equation*}
\int\limits_{\left\vert u\right\vert =1}\left\vert Q\left( \psi \left(
u\right) \right) \right\vert \left\vert du\right\vert \leq \sqrt{%
\int\limits_{E}\frac{\left\vert \varphi ^{\prime }\left( t\right)
\right\vert ^{2}}{\rho \left( t\right) }\left\vert dt\right\vert }\sqrt{%
\int\limits_{E}\left\vert Q\left( t\right) \right\vert ^{2}\rho \left(
t\right) \left\vert dt\right\vert }<\infty ,
\end{equation*}%
which completes the proof.
\end{proof}

\begin{corollary}
For $p,q>1$ such that $\frac{1}{p}+\frac{1}{q}=1$, and a weight function $%
\rho (t)$ on a contour $E,$ let 
\begin{equation}
\delta =\sqrt[q]{\int\limits_{E}\frac{\left\vert \varphi ^{\prime }\left(
t\right) \right\vert ^{\frac{p+q}{p}}}{\rho ^{\frac{q}{p}}\left( t\right) }%
\left\vert dt\right\vert }.  \label{Delta}
\end{equation}%
If $\delta <\infty $, then for any function $Q\left( z\right) \in
L^{p}\left( E,\rho \right) ,$ 
\begin{equation}
\int\limits_{\left\vert u\right\vert =1}\left\vert Q\left( \psi \left(
u\right) \right) \right\vert \left\vert du\right\vert \leq \delta \left\Vert
Q\right\Vert _{L^{p}\left( E,\rho \right) }  \label{Inequtwo}
\end{equation}%
where%
\begin{equation*}
\left\Vert Q\right\Vert _{L^{p}\left( E,\rho \right) }=\sqrt{%
\int\limits_{E}\left\vert Q\left( t\right) \right\vert ^{p}\rho \left(
t\right) \left\vert dt\right\vert }
\end{equation*}
\end{corollary}

\begin{proof}
We have 
\begin{equation*}
\int\limits_{\left\vert u\right\vert =1}\left\vert Q\left( \psi \left(
u\right) \right) \right\vert \left\vert du\right\vert
=\int\limits_{\left\vert u\right\vert =1}\left\vert Q\left( \psi \left(
u\right) \right) \right\vert \frac{\sqrt[p]{\rho \left( \psi \left( u\right)
\right) \left\vert \psi ^{\prime }\left( u\right) \right\vert }}{\sqrt[p]{%
\rho \left( \psi \left( u\right) \right) \left\vert \psi ^{\prime }\left(
u\right) \right\vert }}\left\vert du\right\vert
\end{equation*}

\begin{equation*}
\leq \sqrt[q]{\int\limits_{\left\vert u\right\vert =1}\frac{\left\vert
du\right\vert }{\rho ^{qp^{-1}}\left( \psi \left( u\right) \right)
\left\vert \psi ^{\prime }\left( u\right) \right\vert ^{qp^{-1}}}}\sqrt[p]{%
\int\limits_{\left\vert u\right\vert =1}\left\vert Q\left( \psi \left(
u\right) \right) \right\vert ^{p}\rho \left( \psi \left( u\right) \right)
\left\vert \psi ^{\prime }\left( u\right) \right\vert \left\vert
du\right\vert }
\end{equation*}%
since $\dfrac{1}{p}+\dfrac{1}{q}=1.$ This inequality may be written in the
form%
\begin{equation*}
\int\limits_{\left\vert u\right\vert =1}\left\vert Q\left( \psi \left(
u\right) \right) \right\vert \left\vert du\right\vert \leq \sqrt[q]{%
\int\limits_{E}\frac{\left\vert \varphi ^{\prime }\left( t\right)
\right\vert ^{1+qp^{-1}}}{\rho ^{qp^{-1}}\left( t\right) }\left\vert
dt\right\vert }\sqrt[p]{\int\limits_{E}\left\vert Q\left( t\right)
\right\vert ^{p}\rho \left( t\right) \left\vert dt\right\vert }.
\end{equation*}

This completes the proof.
\end{proof}

\begin{corollary}
If $1\leq q,q,r<\infty $,such that 
\begin{equation*}
\dfrac{1}{q}+\dfrac{1}{p}=\dfrac{1}{r}
\end{equation*}%
For any function $Q\left( z\right) \in L^{p}\left( E,\rho \right) $%
\begin{equation}
\sqrt[r]{\dint\limits_{\left\vert u\right\vert =1}\left\vert Q\left( \psi
\left( u\right) \right) \right\vert ^{r}\left\vert du\right\vert }\leq
\left\Vert Q\right\Vert _{L^{p}\left( E,\rho \right) }\sqrt[q]{%
\dint\limits_{E}\frac{\left\vert \varphi ^{\prime }\left( t\right)
\right\vert ^{1+\tfrac{q}{p}}}{\rho ^{\tfrac{q}{p}}\left( t\right) }%
\left\vert dt\right\vert }  \label{HI5}
\end{equation}%
where%
\begin{equation*}
\left\Vert Q\right\Vert _{L^{p}\left( E,\rho \right) }=\sqrt{%
\dint\limits_{E}\left\vert Q\left( t\right) \right\vert ^{p}\rho \left(
t\right) \left\vert dt\right\vert }
\end{equation*}
\end{corollary}

\begin{proof}
Since the function $Q\left( z\right) $ is defined on the contour $E$, the
function $Q\left( \psi \left( u\right) \right) $ is defined on the
circonference $\left\vert u\right\vert =1.$Because $\left( \dfrac{q}{r}%
\right) ^{-1}+\left( \dfrac{p}{r}\right) ^{-1}=1,$we obtain%
\begin{equation*}
\sqrt[r]{\dint\limits_{\left\vert u\right\vert =1}\left\vert Q\left( \psi
\left( u\right) \right) \right\vert ^{r}\left\vert du\right\vert }=
\end{equation*}

\begin{equation*}
\sqrt[r]{\dint\limits_{\left\vert u\right\vert =1}\left\vert Q\left( \psi
\left( u\right) \right) \right\vert ^{r}\frac{\sqrt[p]{\rho ^{r}\left( \psi
\left( u\right) \right) \left\vert \psi ^{\prime }\left( u\right)
\right\vert ^{r}}}{\sqrt[p]{\rho ^{r}\left( \psi \left( u\right) \right)
\left\vert \psi ^{\prime }\left( u\right) \right\vert ^{r}}}\left\vert
du\right\vert }
\end{equation*}%
\begin{equation*}
\leq \sqrt[q]{\dint\limits_{\left\vert u\right\vert =1}\rho ^{-\tfrac{qr}{pr}%
}\left( \psi \left( u\right) \right) \left\vert \psi ^{\prime }\left(
u\right) \right\vert ^{-\tfrac{qr}{pr}}\left\vert du\right\vert }\sqrt[p]{%
\dint\limits_{\left\vert u\right\vert =1}\left\vert Q\left( \psi \left(
u\right) \right) \right\vert ^{\tfrac{rp}{r}}\rho \left( \psi \left(
u\right) \right) \left\vert \psi ^{\prime }\left( u\right) \right\vert
\left\vert du\right\vert }
\end{equation*}%
Hence%
\begin{equation*}
\sqrt[r]{\dint\limits_{\left\vert u\right\vert =1}\left\vert Q\left( \psi
\left( u\right) \right) \right\vert ^{r}\left\vert du\right\vert }
\end{equation*}%
\begin{equation*}
\leq \sqrt[q]{\dint\limits_{\left\vert u\right\vert =1}\rho ^{-\tfrac{q}{p}%
}\left( \psi \left( u\right) \right) \left\vert \psi ^{\prime }\left(
u\right) \right\vert ^{-\tfrac{q}{p}}\left\vert du\right\vert }\sqrt[p]{%
\dint\limits_{\left\vert u\right\vert =1}\left\vert Q\left( \psi \left(
u\right) \right) \right\vert ^{p}\rho \left( \psi \left( u\right) \right)
\left\vert \psi ^{\prime }\left( u\right) \right\vert \left\vert
du\right\vert }
\end{equation*}%
but,the conformal mapping $\varphi $ of $G$ onto the unit disk $D=\left\{
u\colon \,\left\vert u\right\vert =1\right\} $ such that $\varphi \left( \xi
\right) =0,\varphi ^{\prime }\left( \xi \right) =1$, with $\xi \in G$ fixed,
can be conformally continued by {Caratheodory's} Theorem (see \cite{KO}). If
we suppose that $\psi \colon D=\left\{ u\colon \,\left\vert u\right\vert
=1\right\} \rightarrow G$ is its inverse map, then we have 
\begin{equation*}
\varphi \left( t\right) =u\Longleftrightarrow t=\psi \left( u\right) \text{
\ \ \ \ ,\ \ \ \ \ \ \ }t\in E
\end{equation*}%
and%
\begin{equation*}
\left\vert du\right\vert =\left\vert \varphi ^{\prime }\left( t\right)
\right\vert \left\vert dt\right\vert \text{\ \ \ \ ,\ \ \ \ \ \ \ }t\in E
\end{equation*}%
and%
\begin{equation*}
\varphi ^{\prime }\left( t\right) =\frac{1}{\psi ^{\prime }\left( \varphi
\left( t\right) \right) }\text{\ \ \ \ ,\ \ \ \ \ \ \ }t\in E
\end{equation*}

hence%
\begin{equation*}
\sqrt[r]{\dint\limits_{\left\vert u\right\vert =1}\left\vert Q\left( \psi
\left( u\right) \right) \right\vert ^{r}\left\vert du\right\vert }
\end{equation*}%
\begin{equation*}
\leq \sqrt[q]{\dint\limits_{E}\rho ^{-\tfrac{q}{p}}\left( t\right)
\left\vert \varphi ^{\prime }\left( t\right) \right\vert ^{1+\tfrac{q}{p}%
}\left\vert dt\right\vert }\sqrt{\dint\limits_{E}\left\vert Q\left( t\right)
\right\vert ^{p}\rho \left( t\right) \left\vert dt\right\vert }
\end{equation*}%
The proof of the corollary is completed.
\end{proof}

\end{document}